\documentclass{ifacconf}

\usepackage{graphicx}      
\usepackage{natbib}        
\usepackage{amsmath}
\usepackage{amssymb}
\usepackage{amsfonts}
\usepackage{color}

\usepackage{algorithm}
\usepackage{algorithmicx}
\usepackage{algpseudocode}

\newtheorem{theorem}{Theorem}[section]

\newtheorem{assumption}[theorem]{Assumption}

\begin{document}
\begin{frontmatter}

\title{
    Riccati Recursion for Optimal Control Problems of Nonlinear Switched Systems\thanksref{footnoteinfo}}

\thanks[footnoteinfo]{© 2021 the authors. This work has been accepted to IFAC for publication under a Creative Commons Licence CC-BY-NC-ND”}

\author[First]{S. Katayama} 
\author[First]{T. Ohtsuka} 

\address[First]{
    Department of Systems Science, Graduate School of Informatics, 
    Kyoto University, Sakyo-ku, Kyoto 606-8501, Japan}

\begin{abstract}                
    We propose an efficient algorithm for the optimal control problems (OCPs) of nonlinear switched systems that optimizes the control input and switching instants simultaneously for a given switching sequence.
    We consider the switching instants as the optimization variables and formulate the OCP based on the direct multiple shooting method. 
    We derive a linear equation to be solved in Newton's method and propose a Riccati recursion algorithm to solve the linear equation efficiently.
    The computational time of the proposed method scales linearly with respect to the number of time stages of the horizon as the standard Riccati recursion.
    Numerical experiments show that the proposed method converges with a significantly shorter computational time than the conventional methods.
\end{abstract}


\end{frontmatter}

\section{Introduction}
Switched systems are a class of hybrid systems made up of several dynamical subsystems and the switching laws of active subsystems.
Many practical control systems are modeled as switched systems, such as automobiles with different gears (\cite{bib:automobieGearShift}), electrical circuit systems (\cite{bib:converterMPC}), and mechanical systems with contacts (\cite{bib:bipedWalking, bib:DDP:jumpRobot}).

Optimal control plays a significant role in the planning (e.g., trajectory optimization in robotics) and control, that is, model predictive control (MPC), of dynamical systems, including switched systems. 
There is sufficient literature on the optimal control of linear switched systems, for example, on dynamic programming for linear hybrid systems (\cite{bib:DPForLinearHybridSystem}), mixed integer linear programming (\cite{bib:MIPHybrid}), multiparametric programming (\cite{bib:MultiParametricToolbox}),  or embedding transformation (\cite{bib:SwitchLQRByEmbeddding}).
However, it is generally difficult to solve the optimal control problems (OCPs) for nonlinear switched systems because of a need to solve nonlinear and combinatorial optimization problems. 
A practical research topic of the OCPs of nonlinear switched systems is to optimize the continuous control input and switching instants for a given switching sequence, which reduces the OCP to a continuous optimization problem.
However, it is still difficult to solve such problems within a short computational time, which is crucial for MPC.
\cite{bib:Xu:2004} proposed a two-stage framework in which the upper stage solves the optimization problem to find the optimal switching times under the fixed control input (the solution to the lower stage problem), and the lower stage solves the OCP (both continuous-time OCPs and discrete-time OCPs can be incorporated) under fixed switching times (the solution to the upper stage problem) to find the optimal control input.
However, it was reported that their method required several minutes to converge because they solved a particular form of continuous-time OCP in the lower stage.
\cite{bib:slqSwitch:2017} improved the approach proposed by \cite{bib:Xu:2004}, by solving the discrete-time OCP in the lower stage using a Newton-type method. 
However, their method required several seconds because the problem was still decomposed into two stages as \cite{bib:Xu:2004}, and therefore the real-time application of their algorithm is difficult.

By contrast, our previous study (\cite{bib:Katayama:2020}) applied a Newton-type method to simultaneously optimize the control input and switching time.
It used the Newton-Krylov method and succeeded in the MPC of a compass-like walking robot modeled as a nonlinear switched system with state jumps in real time.
However, the Newton-Krylov method generally requires a careful tuning of settings such as the number of Krylov iterations and preconditioning; otherwise, it can lack numerical stability compared with the direct methods that compute the inverse matrix of the Hessian explicitly. 

In this paper, we propose an efficient algorithm for the OCP of nonlinear switched systems that optimizes the control input and switching instants simultaneously for a given switching sequence.
We formulate the OCP based on the direct multiple shooting method (\cite{bib:directMultipleShooting, bib:realTimeIterationScheme}) with regard to the switching instance as the optimization variable, as well as the state and control input.
We derive the linear equation to be solved in Newton's method and propose a Riccati recursion algorithm to solve the linear equation efficiently.
The computational time of the proposed method scales linearly with respect to the number of time stages of the horizon as the standard Riccati recursion (\cite{bib:Giaf:2016, bib:Nielsen:2017}).
Through numerical experiments, we show that the proposed method converges within a significantly short computational time (10 -- 100 ms), whereas the methods proposed by \cite{bib:Xu:2004} and \cite{bib:slqSwitch:2017} consumed more than several seconds per upper stage iteration.

This paper is organized as follows. 
In Section 2, we formulate the OCP for a nonlinear switched system based on the direct multiple-shooting method. 
In Section 3, we derive the linear equation for Newton’s method and the Riccati recursion algorithm for the linear equation.
In Section 4, we discuss the numerical experiments and demonstrate the effectiveness of the proposed method.
Finally, in Section 5, we conclude our study with a brief summary and a discussion on future work.

\textit{Notation: } We describe the Jacobians and the Hessians of a differentiable function by certain vectors as follows: $\nabla_x f (x)$ denotes $\left( \frac{\partial f}{\partial x} \right) ^{\rm T} (x)$, and $\nabla_{x y} g(x, y)$ denotes $\frac{\partial^2 g}{\partial x \partial y} (x, y)$.

\section{Optimal Control Problems of Switched Systems}\label{section:ocp}
\subsection{Continuous-Time Optimal Control Problem}
We consider the OCP for the following nonlinear switched system consisting of $M$ subsystems
\begin{equation}\label{eq:stateEquation}
    \dot{x} (t) = f_{q(t)} (x(t), u(t)), \;\; f_q : \mathbb{R}^{n_x} \times \mathbb{R}^{n_u} \to \mathbb{R}^{n_x}, 
\end{equation}
where $x (t) \in \mathbb{R}^{n_x}$ denotes the continuous state, $u(t) \in \mathbb{R}^{n_u}$ denotes the piecewise continuous control input, and $q \in \left\{ 1, ..., M \right\}$ denotes the index of the active subsystem.
We also define the switching sequence $\sigma = (q_1, q_2, ..., q_m)$ and the switching time sequence $t_{\sigma} = (t_{1}, t_{2}, ..., t_{m-1})$, where $t_{i}$ denotes an instance of the switch of the active subsystem from subsystem $i$ to subsystem $i+1$ over the time horizon $[t_0, t_f]$.
Note that we do not consider the state jumps at the switch or the state-dependent condition of the switch in this study, the former of which indicates that the state trajectory of (\ref{eq:stateEquation}) is not smooth but continuous over the horizon.
However, in principle, it is possible to extend the proposed method to such cases; the former can be achieved by adding the state jump equation in the proposed Riccati recursion, and the latter can be attained by adding the pure-state constraint representing the switching condition just before the switching instance.
The OCP for the switched systems for a given switching sequence $\sigma$ is defined as follows: find $u$ and the switching time sequence $t_{\sigma}$, that minimize the cost function
\begin{equation}\label{eq:costFunction}
    {J} = \varphi_{q(t_f)} (x(t_f)) + \sum_{i=1}^{m} \int_{t_{i-1}}^{t_{i}} L_{q_i} (x(\tau), u(\tau)) d \tau, 
\end{equation}
where $t_{m} = t_f$,  $\varphi_q (\cdot, \cdot) : \mathbb{R}^{n_x} \to \mathbb{R}$ denotes the terminal cost, and $L_q (\cdot, \cdot)\times \mathbb{R}^{n_x} \times \mathbb{R}^{n_u} \to \mathbb{R}$ denotes the stage cost, subject to (\ref{eq:stateEquation}).

Next, we assume there is only one switch from $q_1$ to $q_2$ over the horizon $[t_0, t_f]$; that is, we assume $\sigma = (q_1, q_2)$ and $t_{\sigma} = (t_s)$, for simplicity.
Note that it is trivial to extend the following formulation and the proposed algorithm to the problems with multiple switches on the horizon (\cite{bib:Katayama:2020}), which we will further discuss in Section 3.

\subsection{Direct Multiple Shooting}
For numerical computation, we discretize the OCP based on the direct multiple shooting method (\cite{bib:directMultipleShooting, bib:realTimeIterationScheme}) and the forward Euler method. 
We divide the horizon $[t_0, t_f]$ into $N$ steps, define the time step $\Delta \tau = (t_f - t_0) / N$, and introduce $i_s$ as an integer that satisfies $i_{s} \Delta \tau \leq t_{s} - t_0 < (i_{s} + 1) \Delta \tau$.
Here, $i_s$ denotes the time stage at which the switch occurs, that is, the switch occurs between stages $i_s$ and $i_{s}+1$.
We introduce the state on the horizon as $x_0$, ..., $x_N$ and control input $u_0$, ..., $u_{N-1}$, where $x_i$ and $u_i$ correspond to $x(t_0 + i \Delta \tau)$ and $u(t_0 + i \Delta \tau)$, respectively.
Note that we include $x_0$ in the optimization variables so that the proposed method can be combined with the real-time algorithm of MPC (\cite{bib:realTimeIterationScheme}).
We also introduce the variables just after the switch as $x_s := x(t_s)$ and $u_s := u(t_s)$.
The cost function (\ref{eq:costFunction}) is then discretized as 
\begin{align}\label{eq:cost:discretized}
    J = & \; \varphi_{q_2} (x_N) + \sum_{i=i_{s}+1}^{N-1} L_{q_2} (x_i, u_i) \Delta \tau \notag \\ 
    & + L_{q_2} (x_s, u_s) (\Delta \tau - \Delta \tau_{s}) + L_{q_{1}} (x_{i_s}, u_{i_s}) \Delta \tau_s \notag, \\ 
    & + \sum_{i=0}^{i_s -1} L_{q_1} (x_i, u_i) \Delta \tau,
\end{align}
where $\Delta \tau_s := t_s - t_0 - i_s \Delta \tau$.
The state equation is discretized as 
\begin{equation}\label{eq:f1}
    x_i + f_{q_1} (x_i, u_i) \Delta \tau - x_{i+1} = 0, \;\; i \in \left\{1, ..., i_s - 1 \right\},
\end{equation}
\begin{equation}\label{eq:fs1}
    x_{i_s} + f_{q_1} (x_{i_s}, u_{i_s}) \Delta \tau_s - x_{s} = 0, 
\end{equation}
\begin{equation}\label{eq:fs2}
    x_s + f_{q_2} (x_s, u_s) (\Delta \tau - \Delta \tau_s) - x_{i_s+1} = 0, 
\end{equation}
and
\begin{equation}\label{eq:f2}
    x_i + f_{q_2} (x_i, u_i) \Delta \tau - x_{i+1} = 0, \;\; i \in \left\{i_s + 1, ..., N - 1 \right\}.
\end{equation}
The constraint on the initial state, 
\begin{equation}\label{eq:x0}
    x_0 - x (t_0) = 0,
\end{equation}
is also imposed because we regard $x_0$ as the optimization variable.
The discretized OCP is now defined as follows: find $x_0$, ..., $x_N$, $x_s$, $u_0$, ..., $u_{N-1}$, $u_s$, and $t_s$ that minimize the cost function (\ref{eq:cost:discretized}) subject to (\ref{eq:f1})--(\ref{eq:x0}).

\subsection{Optimality Conditions}
To derive the optimality conditions, the necessary conditions for optimal control, we define the Hamiltonian $H_{q} (\cdot, \cdot, \cdot) : \mathbb{R}^{n_x} \times \mathbb{R}^{n_u} \times \mathbb{R}^{n_x} \to \mathbb{R}$ for $q \in \left\{ q_1, q_2 \right\}$ as
\begin{equation}
    H_{q} (x, u, \lambda) := L_{q} (x, u) + \lambda ^{\rm T} f_q (x, u).
\end{equation}
The optimality conditions are then derived by the calculus of variations (\cite{bib:Bryson:1975}) as follows:
\begin{align}\label{eq:Hx1}
    \nabla_x H_{q_1} (x_i, u_i, \lambda_{i+1}) \Delta \tau + \lambda_{i+1} - \lambda_i = 0, \notag \\
    i \in \left\{1, ..., i_s - 1 \right\},
\end{align}
\begin{equation}\label{eq:Hxs1}
    \nabla_x H_{q_1} (x_{i_s}, u_{i_s}, \lambda_{s}) \Delta \tau_s + \lambda_{s} - \lambda_{i_s} = 0,
\end{equation}
\begin{equation}\label{eq:Hxs2}
    \nabla_x H_{q_2} (x_{s}, u_{s}, \lambda_{i_s+1}) (\Delta \tau - \Delta \tau_s) + \lambda_{i_s + 1} - \lambda_{s} = 0,
\end{equation}
\begin{align}\label{eq:Hx2}
    \nabla_x H_{q_2} (x_i, u_i, \lambda_{i+1}) \Delta \tau + \lambda_{i+1} - \lambda_i = 0, \notag \\
    i \in \left\{i_s + 1, ..., N - 1 \right\},
\end{align}
\begin{equation}\label{eq:phix}
    \nabla_x \varphi_{q_2} (x_N) - \lambda_N = 0,
\end{equation}
\begin{equation}\label{eq:Hu1}
    \nabla_u H_{q_1} (x_i, u_i, \lambda_{i+1}) \Delta \tau = 0, \;\; i \in \left\{1, ..., i_s - 1 \right\},
\end{equation}
\begin{equation}\label{eq:Hus1}
    \nabla_u H_{q_1} (x_{i_s}, u_{i_s}, \lambda_{s}) \Delta \tau_s = 0,
\end{equation}
\begin{equation}\label{eq:Hus2}
    \nabla_u H_{q_2} (x_{s}, u_{s}, \lambda_{i_s+1}) (\Delta \tau - \Delta \tau_s) = 0,
\end{equation}
\begin{equation}\label{eq:Hu2}
    \nabla_u H_{q_2} (x_i, u_i, \lambda_{i+1}) \Delta \tau = 0, \;\; i \in \left\{i_s + 1, ..., N - 1 \right\},
\end{equation}
and
\begin{equation}\label{eq:H1-H2}
    H_{q_1} (x_{i_s}, u_{i_s}, \lambda_{s}) - H_{q_2} (x_{s}, u_{s}, \lambda_{i_s + 1}) = 0,
\end{equation}
where $\lambda_0, ..., \lambda_N, \lambda_s$ are the Lagrange multipliers with respect to the constraints (\ref{eq:f1})--(\ref{eq:x0}).

\section{Riccati Recursion for Optimal Control Problems of Switched Systems}
\subsection{Linearization for Newton's Method}
The optimality conditions (\ref{eq:f1})--(\ref{eq:H1-H2}) are linearized with $\Delta x_0, ..., \Delta x_N, \Delta x_s$, $\Delta u_0, ..., \Delta u_{N-1}, \Delta u_s$, $\Delta t_{s}$, and $\Delta \lambda_0, \allowbreak ..., \allowbreak \Delta \lambda_N, \Delta \lambda_s$, which are perturbations of $x_0, ..., x_N, x_s$, $u_0, ..., u_{N-1}, u_s$, $t_{s}$, and $\lambda_0, \allowbreak ..., \allowbreak \lambda_N, \lambda_s$, respectively, to apply Newton's method.

\subsubsection{Terminal stage:}
At the terminal stage ($i = N$), we have 
\begin{equation}\label{eq:newton:QxxN}
    Q_{xx, N} \Delta x_N - \Delta \lambda_{N} + \bar{l}_{x, N} = 0,
\end{equation}
where we define $Q_{xx, N} := \nabla_{xx} \varphi_{q_2} (x_N)$ and define $\bar{l}_{x, N}$ using the left-hand side of (\ref{eq:phix}).

\subsubsection{Intermediate stages without a switch:}
At the intermediate stages without a switch, that is, $i < N$ and $i \not= i_s$, we have 
\begin{equation}\label{eq:newton:Qxi}
    Q_{xx, i} \Delta x_i + Q_{xu, i} \Delta u_i + A_i ^{\rm T} \Delta \lambda_{i+1} - \Delta \lambda_{i} + \bar{l}_{x, i} = 0,
\end{equation}
\begin{equation}\label{eq:newton:Qui}
    Q_{xu, i} ^{\rm T} \Delta x_i + Q_{uu, i} \Delta u_i + B_i ^{\rm T} \Delta \lambda_{i+1} + \bar{l}_{u, i} = 0,
\end{equation}
and
\begin{equation}\label{eq:newton:Fi}
    A_i \Delta x_i + B_i \Delta u_i - \Delta x_{i+1} + \bar{x}_i = 0,
\end{equation}
where we define
$Q_{xx, i} := \nabla_{xx} H_q (x_i, u_i, \lambda_{i+1}) \Delta \tau$,
$Q_{xu, i} := \nabla_{xu} H_q (x_i, u_i, \lambda_{i+1}) \Delta \tau$,
$Q_{uu, i} := \nabla_{uu} H_q (x_i, u_i, \lambda_{i+1}) \Delta \tau$,
$A_i := \nabla_x f_q(x_i, u_i) \Delta \tau$,
and
$B_i := \nabla_u f_q(x_i, u_i) \Delta \tau$,
with $q = q_1$ for $i \in \left\{ 1, ..., i_s -1 \right\}$ and $q = q_2$ for $i \in \left\{ i_s + 1 , ..., N -1 \right\}$.
Furthermore, we define $\bar{x}_i$, $\bar{l}_{x, i}$, and $\bar{l}_{u, i}$ using the left-hand sides of (\ref{eq:f1}), (\ref{eq:Hx1}), and (\ref{eq:Hu1}) for $i \in \left\{ 1, ..., i_s -1 \right\}$ and the left-hand sides of (\ref{eq:f2}), (\ref{eq:Hx2}), and (\ref{eq:Hu2}) for $i \in \left\{ i_s + 1 , ..., N -1 \right\}$.

\subsubsection{Intermediate stage with a switch:}
At stage $s$ (switching instant $t_s$), we have 
\begin{align}\label{eq:newton:Qxs}
    Q_{xx, s} \Delta x_s + Q_{xu, s} \Delta u_s + A_s ^{\rm T} \Delta \lambda_{i_s+1} - \Delta \lambda_{s} \notag \\ 
    - h_{x, s} \Delta t_s + \bar{l}_{x, s} = 0,
\end{align}
\begin{equation}\label{eq:newton:Qus}
    Q_{xu, s} ^{\rm T} \Delta x_s + Q_{uu, s} \Delta u_s + B_s ^{\rm T} \Delta \lambda_{i_s+1} - h_{u, s} \Delta t_s + \bar{l}_{u, s} = 0,
\end{equation}
and
\begin{equation}\label{eq:newton:Fs}
    A_s \Delta x_s + B_s \Delta u_s - \Delta x_{i_s+1} - f_s \Delta t_s + \bar{x}_s = 0,
\end{equation}
where we define 
$Q_{xx, s} := \nabla_{xx} H_{q_2} (x_s, u_s, \lambda_{i_s+1}) (\Delta \tau - \Delta \tau_s)$,
$Q_{xu, s} := \nabla_{xu} H_{q_2} (x_s, u_s, \lambda_{i_s+1}) (\Delta \tau - \Delta \tau_s)$,
$Q_{uu, s} := \nabla_{uu} H_{q_2} (x_s, u_s, \lambda_{i_s+1}) (\Delta \tau - \Delta \tau_s)$,
$A_s := \nabla_x f_{q_2} (x_s, u_s) (\Delta \tau \allowbreak - \Delta \tau_s)$,
and 
$B_s := \nabla_u f_{q_2} (x_s, u_s) (\Delta \tau - \Delta \tau_s)$.
We also define $\bar{x}_{i_s}$, $\bar{l}_{x, i_s}$, and $\bar{l}_{u, i_s}$ using the left-hand sides of (\ref{eq:fs2}), (\ref{eq:Hxs2}), and (\ref{eq:Hus2}), and $h_{x, s} := \nabla_x H(x_s, u_s, \lambda_{i_s + 1})$, $h_{u, s} := \nabla_u H(x_s, u_s, \lambda_{i_s + 1})$, and $f_{s} := f_{q_2} (x_s, u_s)$.
At stage $i_s$, we have 
\begin{align}\label{eq:newton:Qxis}
    Q_{xx, i_s} \Delta x_{i_s} + Q_{xu, i_s} \Delta u_{i_s} + A_{i_s} ^{\rm T} \Delta \lambda_{s} - \Delta \lambda_{i_s} \notag \\
    + h_{x, i_s} \Delta t_s + \bar{l}_{x, i_s} = 0,
\end{align}
\begin{equation}\label{eq:newton:Quis}
    Q_{xu, i_s} ^{\rm T} \Delta x_{i_s} + Q_{uu, i_s} \Delta u_{i_s} + B_{i_s} ^{\rm T} \Delta \lambda_{s} + h_{u, i_s} \Delta t_s + \bar{l}_{u, i_s} = 0,
\end{equation}
and
\begin{equation}\label{eq:newton:Fis}
    A_{i_s} \Delta x_{i_s} + B_{i_s} \Delta u_{i_s} - \Delta x_{s} + f_{i_s} \Delta t_s + \bar{x}_{i_s} = 0,
\end{equation}
where we define 
$Q_{xx, i_s} := \nabla_{xx} H_{q_1} (x_{i_s}, u_{i_s}, \lambda_{s}) \Delta \tau_s$,
$Q_{xu, i_s} \allowbreak := \nabla_{xu} H_{q_1} (x_{i_s}, u_{i_s}, \lambda_{s}) \Delta \tau_s$,
$Q_{uu, i_s} := \allowbreak \nabla_{uu} H_{q_1} (x_{i_s}, u_{i_s}, \lambda_{s}) \allowbreak \Delta \tau_s$,
$A_{i_s} := \nabla_x f_{q_1} (x_{i_s}, u_{i_s}) \Delta \tau_s$,
and
$B_{i_s} := \nabla_u f_{q_1} (x_{i_s}, u_{i_s}) \allowbreak \Delta \tau_s$.
We also define $\bar{x}_{i_s}$, $\bar{l}_{x, i_s}$, and $\bar{l}_{u, i_s}$ using the left-hand sides of (\ref{eq:fs1}), (\ref{eq:Hxs1}), and (\ref{eq:Hus1}) and $h_{x, i_s} := \nabla_x H(x_{i_s}, u_{i_s}, \lambda_{s})$, $h_{u, i_s} := \nabla_u H(x_{i_s}, u_{i_s}, \lambda_{s})$, and $f_{i_s} := f_{q_1} (x_{i_s}, u_{i_s})$.
We also have 
\begin{align}\label{eq:newton:h}
    & h_{x, i_s} ^{\rm T} \Delta x_{i_s}
    + h_{u, i_s} ^{\rm T} \Delta u_{i_s}  
    + f_{i_s} ^{\rm T} \Delta \lambda_s \notag \\ 
    & - h_{x, s} ^{\rm T} \Delta x_s  
    - h_{u, s} ^{\rm T} \Delta u_s  
    - f_{s} ^{\rm T} \Delta \lambda_{i_s+1} + \bar{h}_s = 0,
\end{align}
where $\bar{h}_s$ is the left-hand side of (\ref{eq:H1-H2}).

\subsubsection{Initial stage:}
Finally, we have 
\begin{equation}\label{eq:newton:x0}
    \Delta x_0 + \bar{x} = 0,
\end{equation}
where $\bar{x}$ is the left-hand side of (\ref{eq:x0}).

Newton's method for the OCP is reduced to a linear equation to find $\Delta x_0, ..., \Delta x_N, \Delta x_s$, $\Delta u_0, ..., \Delta u_{N-1}, \Delta u_s$, $\Delta t_{s}$, and $\Delta \lambda_0, \allowbreak ..., \allowbreak \Delta \lambda_N, \Delta \lambda_s$, satisfying (\ref{eq:newton:QxxN})--(\ref{eq:newton:x0}).
Next, we make the following reasonable assumption at each Newton iteration:
\assumption{
The linear independence constraint qualification (LICQ) and the second-order sufficient condition (SOSC) hold at each iteration.
}

Note that the assumption on the LICQ and SOSC is equivalent to the positive definiteness of the reduced Hessian (\cite{bib:nocedal}).
To describe this point, let us summarize the stack of the primal Newton directions $\Delta x_0, ..., \Delta x_N, \Delta x_s$, $\Delta u_0, ..., \Delta u_{N-1}, \Delta u_s$, and $\Delta t_{s}$ as $w$ and the stack of the linearized constraints (\ref{eq:newton:Fi}), (\ref{eq:newton:Fs}), (\ref{eq:newton:Fis}), and (\ref{eq:newton:x0}) as $G w + g = 0$.
Then, we describe a quadratic program (QP) that corresponds to the Newton iteration as 
\begin{equation}\label{eq:QP}
    \min_{w} \frac{1}{2} w^{\rm T} H w + l ^{\rm T} w, \;\;\; {\rm s. t.} \;\; G w + g = 0.
\end{equation}
Then, Assumption 3.1 can be considered to be equivalent to the positive definiteness of the reduced Hessian of (\ref{eq:QP}) (\cite{bib:nocedal}).

\subsection{Derivation of Riccati Recursion}
We derive the Riccati recursion to solve the linear equation for Newton's method (\ref{eq:newton:QxxN})--(\ref{eq:newton:x0}).
As the Riccati recursion for the standard OCP (\cite{bib:Giaf:2016, bib:Nielsen:2017}), our goal is the series of matrices $P_i$ and vectors $z_i$ such that $\Delta \lambda_i = P_i \Delta x_i - z_i$ holds.

\subsubsection{Terminal stage:}
At the terminal stage ($i = N$), 
\begin{equation}\label{eq:riccati:N}
    P_N = Q_{xx, N}, \; z_N = - \bar{l}_N
\end{equation}
is given as the standard Riccati recursion. 
In the forward recursion, we have $\Delta x_N$ and compute $\Delta \lambda_N = P_N \Delta x_N - z_N$.

\subsubsection{Intermediate stages without a switch:}
At the intermediate stages without a switch ($i < N$ and $i \not= i_s$), the following standard backward Riccati recursion is given under Assumption 3.1 and an assumption that we have $P_{i+1}$ and $z_{i+1}$ satisfying $\Delta \lambda_{i+1} = P_{i+1} \Delta x_{i+1} - z_{i+1}$ (\cite{bib:Giaf:2016, bib:Nielsen:2017}):
\begin{equation}\label{eq:riccati:F}
    F_i := Q_{xx, i} + A_i ^{\rm T} P_{i+1} A_i,
\end{equation}
\begin{equation}\label{eq:riccati:H}
    H_i := Q_{xu, i} + A_i ^{\rm T} P_{i+1} B_i,
\end{equation}
\begin{equation}\label{eq:riccati:G}
    G_i := Q_{uu, i} + B_i ^{\rm T} P_{i+1} B_i,
\end{equation}
\begin{equation}\label{eq:riccati:K}
    K_i := - G_i ^{-1} H_i ^{\rm T}, \; k_i := - G_i ^{-1} (B_i ^{\rm T} P_{i+1} \bar{x}_i - B_i ^{\rm T} z_{i+1} + \bar{l}_{u, i}),
\end{equation}
and
\begin{equation}\label{eq:riccati:P}
    P_i := F_i - K_i ^{\rm T} G_{i} K_i, \; z_i := A_i ^{\rm T} (z_{i+1} - P_{i+1} \bar{x}_i) - \bar{l}_{x, i} - H_i k_i.
\end{equation}
In the forward recursion, we have $\Delta x_i$ and compute $\Delta u_i$ and $\Delta \lambda_i$ from $\Delta x_i$ as
\begin{equation}\label{eq:riccati:ui}
    \Delta u_i = K_i \Delta x_i + k_i,  
\end{equation}
\begin{equation}\label{eq:riccati:lmdi}
    \Delta \lambda_i = P_i \Delta x_i - z_i,
\end{equation}
and compute $\Delta x_{i+1}$ from (\ref{eq:newton:Fi}).

\subsubsection{Intermediate stages with a switch:}
Next, we derive the Riccati recursion for intermediate stages with a switch, that is, at stages $s$ and $i_s$.
Suppose that we have $P_{i_s+1}$ and $z_{i_s+1}$ satisfying (\ref{eq:riccati:lmdi}).
To obtain the Riccati recursion, we derive the relation of $\Delta \lambda_{i_s}$, $\Delta u_{i_s}$, $\Delta \lambda_s$, $\Delta x_s$, $\Delta u_s$, $\Delta \lambda_{i_s+1}$, and $\Delta t_s$ with respect to $\Delta x_{i_s}$. This problem is equivalent to factorizing the linear equation (\ref{eq:linearProblemSwitch}).
\begin{figure*}
\begin{equation}\label{eq:linearProblemSwitch}
    \begin{bmatrix}
        - I & Q_{xx, i_s} & Q_{xu, i_s} & A_{i_s} ^{\rm T} & h_{x, i_s} \\ 
        & Q_{xu, i_s} ^{\rm T} & Q_{uu, i_s} & B_{i_s} ^{\rm T} & h_{u, i_s} \\ 
        & A_{i_s} & B_{i_s} &   & f_{i_s} & - I \\ 
        & h_{x, i_s} ^{\rm T} & h_{u, i_s} ^{\rm T} & f_{i_s} ^{\rm T} &  & - h_{x, s} ^{\rm T} & - h_{u, s} ^{\rm T} & - f_{s} ^{\rm T} \\ 
        & & & -I & - h_{x, s} & Q_{xx, s} & Q_{xu, s} & A_{s} ^{\rm T} \\ 
        & & &    & - h_{u, s} & Q_{xu, s} ^{\rm T} & Q_{uu, s} & B_{s} ^{\rm T} \\ 
        & & &    & - f_{s}    & A_{s} & B_{s} &   & - I   \\ 
        & & &    &            &   &  & - I & P_{i_s+1} 
    \end{bmatrix}
    \begin{bmatrix}
        \Delta \lambda_{i_s} \\
        \Delta x_{i_s} \\
        \Delta u_{i_s} \\
        \Delta \lambda_s \\
        \Delta t_s \\
        \Delta x_s \\
        \Delta u_s \\
        \Delta \lambda_{i_s+1} \\
        \Delta x_{i_s+1} \\
    \end{bmatrix}
    = - \begin{bmatrix}
        \bar{l}_{x, i_s} \\
        \bar{l}_{u, i_s} \\
        \bar{x}_{i_s} \\
        \bar{h}_s \\
        \bar{l}_{x, s} \\
        \bar{l}_{u, s} \\
        \bar{x}_{s} \\
        - z_{i_s+1} \\
    \end{bmatrix}
\end{equation}
\end{figure*}
First, by using (\ref{eq:newton:Fs}), (\ref{eq:riccati:lmdi}) with $i = i_s+1$, and 
\begin{equation}\label{eq:riccati:us}
    \Delta u_s = K_s \Delta x_s + k_s - T_s \Delta t_s ,
\end{equation}
we can reduce (\ref{eq:linearProblemSwitch}) into 
\begin{equation}\label{eq:linearProblemSwitchReduced1}
    \begin{bmatrix}
        - I & Q_{xx, i_s} & Q_{xu, i_s} & A_{i_s} ^{\rm T} & h_{x, i_s} \\ 
         & Q_{xu, i_s} ^{\rm T} & Q_{uu, i_s} & B_{i_s} ^{\rm T} & h_{u, i_s} \\ 
         & A_{i_s} & B_{i_s} &   & f_{i_s} & - I \\ 
         & h_{x, i_s} ^{\rm T} & h_{u, i_s} ^{\rm T} & f_{i_s} ^{\rm T} & \xi_s & - \Gamma_{s} ^{\rm T} & \\ 
         & & & -I & - \Gamma_{s} & P_{s} 
    \end{bmatrix}
    \begin{bmatrix}
        \Delta \lambda_{i_s} \\
        \Delta x_{i_s} \\ 
        \Delta u_{i_s} \\
        \Delta \lambda_s \\
        \Delta t_s \\
        \Delta x_s 
    \end{bmatrix} 
    = 
    - \begin{bmatrix}
        \bar{l}_{x, {i_s}} \\ 
        \bar{l}_{u, {i_s}} \\ 
        \bar{x}_{{i_s}} \\ 
        \eta_{s} \\ 
        - z_s
    \end{bmatrix} 
\end{equation}
where we define $F_s$, $H_s$, $G_s$, $K_s$, $k_s$, $P_s$, and $z_s$ according to (\ref{eq:riccati:F})--(\ref{eq:riccati:P}) 
and
\begin{equation}\label{eq:riccati:psis}
    \psi_{x, s} := h_{x, s} + A_s ^{\rm T} P_{i_s+1} f_s, \;\;  \psi_{u, s} := h_{u, s} + B_s ^{\rm T} P_{i_s+1} f_s,
\end{equation}
\begin{equation}\label{eq:riccati:Ts}
    T_{s} := - G_s ^{-1} \psi_{u, s}, \;\; \Gamma_{s} := \psi_{x, s} + K_{s} ^{\rm T} \psi_{u, s},
\end{equation}
\begin{equation}\label{eq:riccati:xis}
    \xi_s := f_s ^{\rm T} P_{i_s+1} f_s - \psi_{u, s} ^{\rm T} G_{s} ^{-1} \psi_{u, s},
\end{equation}
and
\begin{equation}\label{eq:riccati:etas}
    \eta_s := \bar{h}_{s} - f_{s} ^{\rm T} (P_{i_s+1} \bar{x}_s - z_{i_s+1}) - \psi_{u, s} k_s .
\end{equation}
We further factorize (\ref{eq:linearProblemSwitchReduced1}) using (\ref{eq:newton:Fis}),
\begin{equation}\label{eq:riccati:lmds}
    \Delta \lambda_s = P_s \Delta x_s - z_s - \Gamma_s \Delta t_s 
\end{equation}
and
\begin{equation}\label{eq:riccati:uis}
    \Delta u_{i_s} = K_{i_s} \Delta x_{i_s} + k_{i_s} + T_{i_s} \Delta t_s ,
\end{equation}
and then obtain
\begin{align}
    \label{eq:linearProblemSwitchReduced2}
    & \begin{bmatrix}
        - I & F_{i_s} - K_{i_s} ^{\rm T} G_{i_s} K_{i_s} & \Gamma_{i_s} \\ 
            & \Gamma_{i_s} ^{\rm T} & \tilde{\xi}_s
    \end{bmatrix}
    \begin{bmatrix}
        \Delta \lambda_{i_s} \\ 
        \Delta x_{i_s} \\ 
        \Delta t_s
    \end{bmatrix} \notag \\
    & = - \begin{bmatrix}
        l_{x, {i_s}} + A_{i_s} ^{\rm T} (P_s \bar{x}_{i_s} - s_s) + H_{i_s} k_{i_s} \\ 
        \tilde{\eta}_{s} 
    \end{bmatrix} ,
\end{align}
where we define $F_{i_s}$, $H_{i_s}$, $G_{i_s}$, $K_{i_s}$, and $k_{i_s}$ according to (\ref{eq:riccati:F})--(\ref{eq:riccati:K}) and
\begin{equation}\label{eq:riccati:psixis}
    \psi_{x, {i_s}} := h_{x, {i_s}} + A_{i_s} ^{\rm T} P_{s} f_{i_s} - A_{i_s} ^{\rm T} \Gamma_s, 
\end{equation}
\begin{equation}\label{eq:riccati:psiuis}
    \psi_{u, {i_s}} := h_{u, {i_s}} + B_{i_s} ^{\rm T} P_{s} f_{i_s} - B_{i_s} ^{\rm T} \Gamma_s,
\end{equation}
\begin{equation}\label{eq:riccati:Tis}
    T_{{i_s}} := - G_{i_s} ^{-1} \psi_{u, {i_s}}, \;\; \Gamma_{{i_s}} := \psi_{x, {i_s}} + K_{{i_s}} ^{\rm T} \psi_{u, {i_s}},
\end{equation}
\begin{equation}\label{eq:riccati:xiis}
    \tilde{\xi}_s = \xi_s - 2 \Gamma_s ^{\rm T} f_{i_s} + f_{i_s} ^{\rm T} P_s f_{i_s} - \psi_{u, k} ^{\rm T} G_{{i_s}} ^{-1} \psi_{u, {i_s}},
\end{equation}
and 
\begin{equation}\label{eq:riccati:etais}
    \tilde{\eta}_s = \eta_s - \Gamma_s ^{\rm T} \bar{x}_{i_s} + f_{i_s} ^{\rm T} (P_s \bar{x}_{i_s} - s_{s} ) + \psi_{x, {i_s}} ^{\rm T} k_{{i_s}}.
\end{equation}
Finally, we have recursions
\begin{equation}\label{eq:riccati:Pis}
    P_{i_s} := F_{i_s} - K_{i_s} ^{\rm T} G_{{i_s}} K_{i_s} - \tilde{\xi}_s ^{-1} \Gamma_{i_s} \Gamma_{i_s} ^{\rm T}
\end{equation}
and
\begin{equation}\label{eq:riccati:zis}
    z_{i_s} := A_{i_s} ^{\rm T} (z_{s} - P_{s} \bar{x}_{i_s}) - \bar{l}_{x, {i_s}} - H_{i_s} k_{i_s} +  \tilde{\xi}_s ^{-1} \Gamma_{i_s} \tilde{\eta}_s.
\end{equation}
We can then compute $\Delta t_s$ by
\begin{equation}\label{eq:riccati:ts}
    \Delta t_s = - \tilde{\xi}_s ^{-1} \Gamma_{i_s} ^{\rm T} \Delta x_{i_s} - \tilde{\xi}_s ^{-1} \tilde{\eta}_s
\end{equation}
and $\Delta \lambda_{i_s}$ from (\ref{eq:riccati:lmdi}) with $i = i_s$.
In the forward recursion, we have $\Delta x_{i_s}$ and compute $\Delta t_s$, $\Delta u_{i_s}$, $\Delta \lambda_{i_s}$, and $\Delta x_{s}$ from (\ref{eq:riccati:ts}), (\ref{eq:riccati:uis}), (\ref{eq:riccati:lmdi}), and (\ref{eq:riccati:xis}).
Subsequently, we compute $\Delta u_{s}$, $\Delta \lambda_{s}$, and $\Delta x_{i_s+1}$ from (\ref{eq:riccati:us}), (\ref{eq:riccati:lmds}), and (\ref{eq:newton:Fs}).

Note that $\tilde{\xi}_s > 0$ holds under Assumption 3.1.
By the contraposition of this property, if $\tilde{\xi}_s \leq 0$, Assumption 3.1 does not hold, that is, the stationary point is not a local minimum (a maximum or a saddle point).
This property is verified as follows:
First, recall that the Newton iteration corresponds to solving the QP (\ref{eq:QP}) and the Riccati recursion corresponds to applying dynamic programming (DP) to the QP (\cite{bib:Giaf:2016, bib:Nielsen:2017}).
By Lemma 16.1 of \cite{bib:nocedal} and Assumption 3.1 (positive definiteness of the reduced Hessian), the QP (\ref{eq:QP}) has a unique solution.
By solving this QP backward in time using DP, we obtain the following QP subproblem: find $\Delta u_{i_s}$, $\Delta t_{s}$, $\Delta x_{s}$, and $\Delta u_{s}$ that minimizes a cost function whose quadratic term is given by
\begin{align}\label{eq:QPsub}
    & \frac{1}{2}
    \begin{bmatrix}
        \Delta x_{i_s} \\ 
        \Delta u_{i_s} \\
        \Delta t_{s} \\
        \Delta x_{s} \\
        \Delta u_{s} \\
    \end{bmatrix}^{\rm T}
    \begin{bmatrix}
        Q_{xx, i_s} & Q_{xu, i_s}          & h_{x, i_s} \\
        Q_{xu, i_s} ^{\rm T} & Q_{uu, i_s} & h_{u, i_s}  \\
        h_{x, i_s} ^{\rm T} &  h_{u, i_s} ^{\rm T} & & - h_{x, s} ^{\rm T} & - h_{u, s} ^{\rm T} \\
        & & - h_{x, s} & Q_{xx, s} & Q_{xu, s} \\
        & & - h_{u, s} & Q_{xu, s} ^{\rm T} & Q_{uu, s} 
    \end{bmatrix} 
    \begin{bmatrix}
        \Delta x_{i_s} \\ 
        \Delta u_{i_s} \\
        \Delta t_{s} \\
        \Delta x_{s} \\
        \Delta u_{s} \\
    \end{bmatrix} \notag \\
    & + \frac{1}{2} \Delta x_{i_s+1} ^{\rm T} P_{i_s + 1} \Delta x_{i_s+1},
\end{align}
where the last term originates from the cost-to-go function of stage $i_s + 1$, 
subject to (\ref{eq:newton:Fis}) and (\ref{eq:newton:Fs}).
This QP subproblem has a unique solution as well as (\ref{eq:QP}).
By eliminating $\Delta u_{i_s}$, $\Delta x_{s}$, $\Delta u_{s}$, and $\Delta x_{i_s+1}$ from the above QP (\ref{eq:QPsub}), we obtain another QP: find $\Delta t_s$ that minimizes a cost function whose quadratic term is given by 
\begin{equation*}
\frac{1}{2}
\begin{bmatrix}
    \Delta x_{i_s} \\ 
    \Delta t_{s} 
\end{bmatrix}^{\rm T}
\begin{bmatrix}
    F_{i_s} - K_{i_s} ^{\rm T} G_{i_s} K_{i_s} & \Gamma_{i_s} \\ 
    \Gamma_{i_s} ^{\rm T} & \tilde{\xi}_s
\end{bmatrix}
\begin{bmatrix}
    \Delta x_{i_s} \\ 
    \Delta t_{s}
\end{bmatrix}.
\end{equation*}
Since this QP must have a solution, $\tilde{\xi}_s > 0$ holds.

\subsubsection{Initial stage:}
At the beginning of the forward recursion, we compute $\Delta x_0$ from (\ref{eq:newton:x0}).

\subsection{Step Size Selection}\label{subsec:stepSize}
The full-step Newton's method can be very aggressive at the beginning of Newton's iterations.
In such cases, the magnitude of the switching time direction, $\Delta t_s$, is excessively large, and the solution diverges or converges to the saddle points.
To avoid such situations, we selected the step size based on the switching time.
We assumed that each switching time, $t_s$, must lie on $[t_{s, \min}, t_{s, \max}]$.
Subsequently, we chose the step size $\alpha$ by applying the fraction-to-boundary rule (\cite{bib:Ipopt}) for the inequality constraints $t_s - t_{s, \min} > 0$ and $t_{s, \max} - t_s > 0$.

\subsection{Algorithm and Convergence}

\begin{algorithm}[tb]
\caption{Computation of Newton direction by the proposed Riccati recursion}
\label{alg1}
\begin{algorithmic}[1]
    \Require Initial state ${x} (t_0)$ and the current solution $x_0$, ..., $x_{N}$, $x_s$, $u_0$, ..., $u_{N-1}$, $u_s$, $\lambda_0$, ..., $\lambda_N$, $\lambda_s$, and $t_s$.
    \Ensure Newton directions $\Delta x_0$, ..., $\Delta x_{N}$, $\Delta x_s$, $\Delta u_0$, ..., $\Delta u_{N-1}$, $\Delta u_s$, $\Delta \lambda_0$, ..., $\Delta \lambda_N$, $\Delta \lambda_s$, and $\Delta t_s$.
    \State Form the linear equations, i.e., compute the coefficient matrices and residuals of (\ref{eq:newton:QxxN})--(\ref{eq:newton:x0}).
    \State Compute $P_N$ and $z_N$ from (\ref{eq:riccati:N}). 
    \For{$i=N,\cdots,i_s+1$} 
        \State Compute $P_i$ and $z_i$ from (\ref{eq:newton:Fi})--(\ref{eq:riccati:P}). 
    \EndFor
    \State Compute $P_{s}$, $z_s$, and $P_{i_s}$ and $z_{i_s}$ from (\ref{eq:newton:Fi})--(\ref{eq:riccati:P}), (\ref{eq:riccati:psis})--(\ref{eq:riccati:etas}), and (\ref{eq:riccati:psixis})--(\ref{eq:riccati:zis}), respectively.
    \For{$i=i_s-1,\cdots,0$} 
        \State Compute $P_i$ and $z_i$ from (\ref{eq:newton:Fi})--(\ref{eq:riccati:P}). 
    \EndFor
    \State Compute $x_0$ from (\ref{eq:newton:x0}).
    \For{$i=0,\cdots,i_s-1$} 
        \State Compute $\Delta u_i$, $\Delta \lambda_i$, and $\Delta x_{i+1}$ from (\ref{eq:riccati:ui}), (\ref{eq:riccati:lmdi}), and (\ref{eq:newton:Fi}), respectively.
    \EndFor
    \State Compute $\Delta t_s$, $\Delta u_{i_s}$, $\Delta \lambda_{i_s}$, $\Delta x_{s}$, $\Delta u_{s}$, $\Delta \lambda_{s}$, and $\Delta x_{i_s+1}$, from (\ref{eq:riccati:ts}), (\ref{eq:riccati:uis}), (\ref{eq:riccati:lmdi}), (\ref{eq:newton:Fis}), (\ref{eq:riccati:us}), (\ref{eq:riccati:lmds}), and (\ref{eq:newton:Fs}), respectively.
    \For{$i=i_s+1,\cdots,N-1$} 
        \State Compute $\Delta u_i$, $\Delta \lambda_i$, and $\Delta x_{i+1}$ from (\ref{eq:riccati:ui}), (\ref{eq:riccati:lmdi}), and (\ref{eq:newton:Fi}).
    \EndFor
    \State Compute $\Delta \lambda_N$ from (\ref{eq:riccati:lmdi}).
\end{algorithmic}
\end{algorithm}

We summarize the single Newton iteration using the proposed Riccati recursion algorithm in Algorithm \ref{alg1}.
As shown in Algorithm \ref{alg1}, the proposed Riccati recursion computes the Newton direction for a given solution. 
In the first step, we form the linear equations of Newton's method, that is, compute the coefficient matrices and residuals of (\ref{eq:newton:QxxN})--(\ref{eq:newton:x0}) (line 1).
Second, we perform the backward Riccati recursion and compute $P_i$ and $z_i$ for $i \in \left\{ 1, ..., N, s \right\}$ (lines 3--9).
Finally, we perform the forward Riccati recursion and compute the Newton directions for all the variables (lines 10--18).

We also summarize Newton's method for the OCP of the switched system in Algorithm \ref{alg2}.
We iterate Algorithm \ref{alg1} and update the solution until the norm of residuals of the optimality conditions (\ref{eq:f1})--(\ref{eq:H1-H2}) become smaller than a prespecified threshold ($\epsilon$ in Algorithm \ref{alg2}).
We refer to the $l_2$-norm of the residual of the optimality conditions as ``Opt. error'' in the following sections.
After each iteration of the proposed Riccati recursion (line 3), we determine the step size (line 4), for example, by the fraction-to-boundary rule introduced in subsection \ref{subsec:stepSize}, and update all variables (line 5).
After updating the solution, we update $i_s$ based on $t_s ^{k+1}$ (line 6).
Therefore, $i_s$ at iteration $k+1$ can vary from that at iteration $k$.

Note that the proposed method can be directly applied to the OCP with multiple switches on the horizon.
When there are multiple switches on the horizon, we compute the coefficient matrices and residuals in (\ref{eq:newton:Qxs})--(\ref{eq:newton:h}) for each switch in line 1 of Algorithm \ref{alg1}, apply line 6 of Algorithm \ref{alg1} for each switch in the backward Riccari recursion, and apply line 14 of Algorithm \ref{alg1} for each switch in the forward Riccari recursion.

The proposed method is in substance a Newton's method for the optimization problem with equality constraints as the standard direct multiple shooting method, and we can similarly discuss the convergence.
The difference is that the switching stage, $i_s$, can change depending on $t_s$ at each iteration, which implies that the problem structure, that is, the cost function (\ref{eq:costFunction}) and equality constraints  (\ref{eq:f1})--(\ref{eq:f2}), can change between the iterations.
If $i_s$ does not change, that is, the optimal switching instant $t_s$ satisfies $i_s \Delta \tau < t_s - t_0 < (i_s + 1) \Delta \tau$, and the sequence of the switching instant at each iteration $\left\{t_s ^k \right\}$ satisfies $i_s \Delta \tau < t_s ^k - t_0 < (i_s + 1) \Delta \tau$ for the same $i_s$ with $\alpha = 1$, then the proposed method achieves quadratic convergence under Assumption 3.1, for example, by Theorem 18.4 of \cite{bib:nocedal}.

\begin{algorithm}[tb]
\caption{Newton's method for the OCP of the switched system by the proposed Riccati recursion}
\label{alg2}
\begin{algorithmic}[1]
    \Require Initial state ${x} (t_0)$, the initial guess of the solution $x_0 ^0, ..., x_{N} ^0, x_s ^0$, $u_0 ^0, ..., u_{N-1} ^0, x_s^0$, $\lambda_0 ^0, ..., \lambda_N ^0, \lambda_s ^0$, and $t_s ^0$, and the termination criteria $\epsilon \geq 0$.
    \Ensure Optimal solution $x_0$, ..., $x_{N}$, $x_s$, $u_0$, ..., $u_{N-1}$, $u_s$, $\lambda_0$, ..., $\lambda_N$, $\lambda_s$, and $t_s$
    \State Set $i_s$ such that it satisfies $i_{s} \Delta \tau \leq t_{s} ^{0} - t_0 < (i_{s} + 1) \Delta \tau$.
    \While{${\rm Opt. error} > \epsilon$} (at $k$-th iteration)
        \State Compute the Newton directions $\Delta x_0 ^k$, ..., $\Delta x_{N} ^k$, $\Delta x_s ^k$, $\Delta u_0 ^k$, ..., $\Delta u_{N-1} ^k$, $\Delta u_s ^k$, $\Delta \lambda_0 ^k$, ..., $\Delta \lambda_N ^k$, $\Delta \lambda_s ^k$, and $\Delta t_s ^k$ using {\bf Algorithm 1} based on the current solution.
        \State Choose step size $\alpha$ ($0 < \alpha \leq 1$).
        \State Update the solution by 
        $x_i ^{k+1} \leftarrow x_i ^k + \alpha \Delta x_i ^k$, 
        $u_i ^{k+1} \leftarrow u_i ^k + \alpha \Delta u_i ^k$, 
        $\lambda_i ^{k+1} \leftarrow \lambda_i ^k + \alpha \Delta \lambda_i ^k$,  and $t_s ^{k+1} \leftarrow t_s ^k + \alpha \Delta t_s ^k$.
        \State Update $i_s$ such that it satisfies $i_{s} \Delta \tau \leq t_{s} ^{k+1} - t_0 < (i_{s} + 1) \Delta \tau$.
    \EndWhile
\end{algorithmic}
\end{algorithm}

\section{Numerical Experiments}
We conducted two numerical experiments to show the effectiveness of the proposed method.
The proposed algorithm was written using Julia, and all the experiments were conducted on a laptop with a quad-core CPU Intel Core i7-10510U @1.8 GHz.

\subsection{Example 1}
The first example is a switched system consisting of two linear subsystems that were treated in \cite{bib:Xu:2004}.
The dynamics of the subsystems are given by 
\begin{equation*}
    f_1 (x, u) = \begin{bmatrix}
          0.6 & 1.2 \\
        - 0.8 & 3.4 
    \end{bmatrix}
    \begin{bmatrix}
     x_1 \\ x_2   
    \end{bmatrix}
    + \begin{bmatrix}
        1 \\
        1
    \end{bmatrix}
    u
\end{equation*}
and
\begin{equation*}
    f_2 (x, u) = \begin{bmatrix}
          4 & 3 \\
        - 1 & 0
    \end{bmatrix}
    \begin{bmatrix}
     x_1 \\ x_2   
    \end{bmatrix}
    + \begin{bmatrix}
        2 \\
        - 1
    \end{bmatrix}
    u,
\end{equation*}
the terminal cost is given by
\begin{equation*}
    \varphi_{q} (x) = \frac{1}{2} (x_1 - x_{\rm ref, 1}  )^2 + \frac{1}{2} (x_2 - x_{\rm ref, 2}  )^2 , \;\; q \in \left\{ 1, 2 \right\},
\end{equation*}
and the stage cost is given by
\begin{equation*}
    L_{q} (x, u) = \frac{1}{2} (x_2 - x_{\rm ref, 2} ) ^2 + \frac{1}{2} u ^2 , \;\; q \in \left\{ 1, 2 \right\},
\end{equation*}
where $t_0 = 0$, $t_f = 2$, $x_{\rm ref, 1} = 4$, and $x_{\rm ref, 2} = 2$.
We discretize the continuous-time OCP into $N = 175$ steps.
The switching sequence is given by $\sigma = (1, 2)$.
The initial state is given by $x(t_0) = [2, 3] ^{\rm T}$, and we set the initial guess of the solution as $x_i = x(t_0)$ for $i \in \left\{ 1, ..., N, s \right\}$ and $t_1 = 1.0$ [s].
In the step-size selection, we use the fraction-to-boundary rule introduced in subsection \ref{subsec:stepSize} such that $t_1 \in [t_0, t_f]$.

\begin{figure}[tb]
\begin{center}
\includegraphics[scale=0.525]{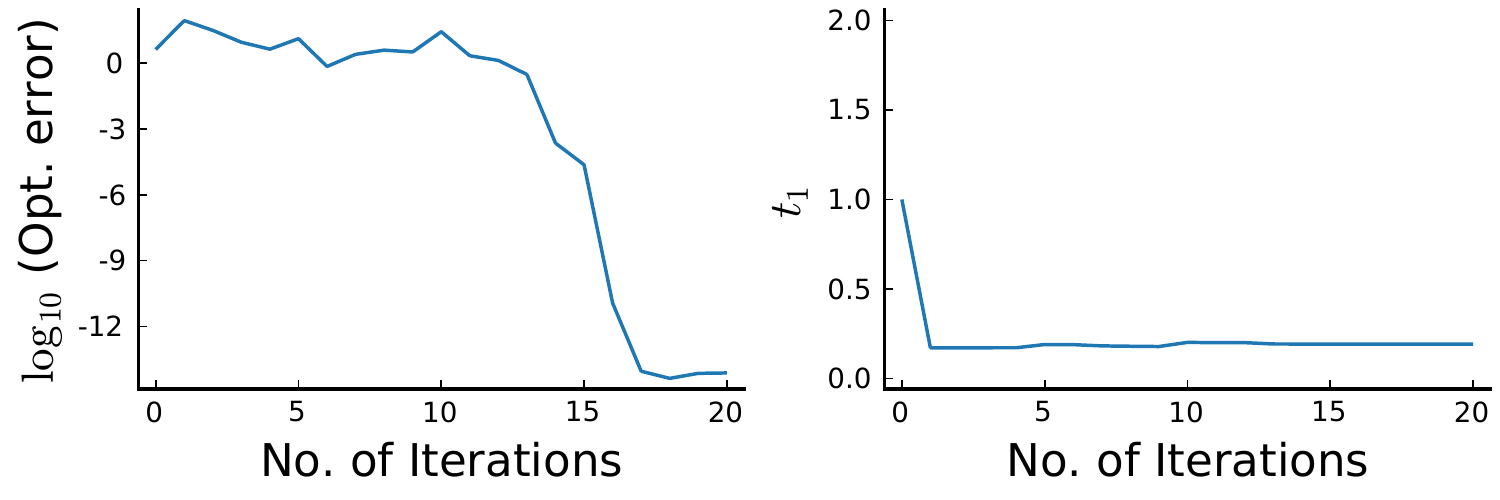}  
\caption{$\log_{10}$ scaled Opt. error and switching instant $t_1$ with respect to the iterations of Example 1.} 
\label{fig:ex1:convergence}
\end{center}
\end{figure}

Figure \ref{fig:ex1:convergence} shows the $\log_{10}$ scaled Opt. error ($l_2$-norm of the residual in the optimality conditions (\ref{eq:f1})--(\ref{eq:H1-H2})), and the switching instants with respect to the number of Newton iterations.
As shown in Fig. \ref{fig:ex1:convergence}, the proposed method converges after 17 Newton iterations.
The two-stage method proposed by \cite{bib:Xu:2004} also converged after 17 upper stage gradient iterations.
However, it took a significantly long computational time (over 30 s using MATLAB on AMD Athlon @4990-MHz) because it solves the continuous-time OCP in the lower stage problem after each upper stage gradient iteration for the switching instants.
By contrast, the proposed method required only 1 ms per Newton iteration and 17 ms for the convergence owing to the discretized formulation and the proposed Riccati recursion algorithm.
The solution provided by the proposed method is almost identical to the continuous-time counterpart reported in \cite{bib:Xu:2004}; for example, the optimal switching instant of the proposed method is 0.1919 [s] and that of \cite{bib:Xu:2004} was 0.1897 [s]. 
In addition, almost the same state and control input trajectory was obtained as \cite{bib:Xu:2004}.

\subsection{Example 2}
The second example is a switched system consisting of three nonlinear subsystems that were treated in \cite{bib:Xu:2004} and \cite{bib:slqSwitch:2017}.
The dynamics of the subsystems are given by 
\begin{equation*}
    f_1 (x, u) = 
    \begin{bmatrix}
    x_1 + u_1 \sin(x_1) \\
    - x_2 - u_1 \cos(x_2) 
    \end{bmatrix}, 
\end{equation*}
\begin{equation*}
    f_2 (x, u) = 
    \begin{bmatrix}
    x_2 + u_1 \sin(x_2) \\
    - x_1 - u_1 \cos(x_1) 
    \end{bmatrix}, 
\end{equation*}
and
\begin{equation*}
    f_3 (x, u) = 
    \begin{bmatrix}
    - x_1 - u_1 \sin(x_1) \\
    x_2 + u_1 \cos(x_2) 
    \end{bmatrix},
\end{equation*}
the terminal cost is given by
\begin{equation*}
    \varphi_{q} (x) = \frac{1}{2} ||x - x_{\rm ref} ||^2 , \;\; q \in \left\{ q_1, q_2, q_3 \right\},
\end{equation*}
and the stage cost is given by
\begin{equation*}
    L_{q} (x, u) = \frac{1}{2} ||x - x_{\rm ref} ||^2 + || u ||^2 , \;\; q \in \left\{ q_1, q_2, q_3 \right\},
\end{equation*}
where $t_0 = 0$, $t_f = 3$, and $x_{\rm ref} = [1, -1]^{\rm T}$.
We discretize the continuous-time OCP into $N=220$ steps.
The switching sequence is given by $\sigma = (1, 2, 3)$.
The initial state is given by $x(t_0) = [2, 3] ^{\rm T}$, and the initial guess of the solution is given by $x_i = x(t_0)$, $[t_{1, 2}, t_{2,3}] = [0.5, 1.0]$.
As in the previous example, we use the fraction-to-boundary rule in subsection \ref{subsec:stepSize} in the step-size selection such that $t_1, t_2 \in [t_0, t_f]$.

\begin{figure}[tb]
\begin{center}
\includegraphics[scale=0.525]{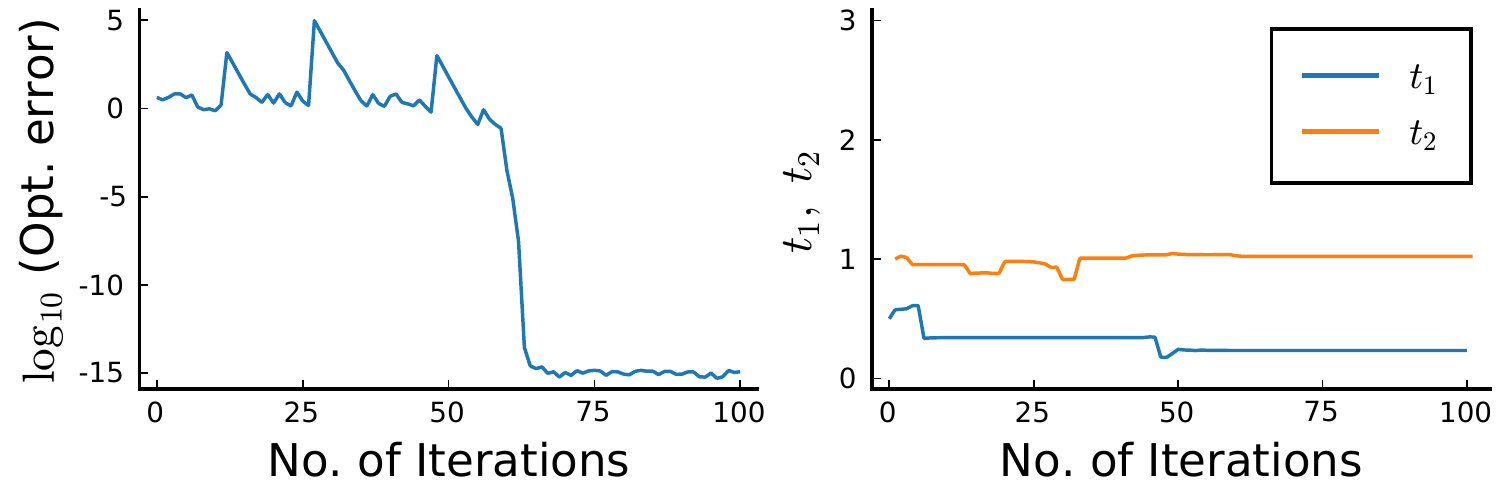}   
\caption{$\log_{10}$ scaled Opt. error and the switching instants $t_1$ and $t_2$ with respect to the iterations of Example 2.} 
\label{fig:ex2:convergence}
\end{center}
\end{figure}

Figure \ref{fig:ex2:convergence} shows the $\log_{10}$ scaled Opt. error and the switching instants with respect to the number of Newton iterations.
As shown in Fig. \ref{fig:ex2:convergence}, the proposed method converges after 70 Newton iterations.
The proposed method takes only 1.153 ms per Newton iteration and 80 ms for convergence, which is significantly faster than the existing methods (\cite{bib:Xu:2004} and \cite{bib:slqSwitch:2017})). 
\cite{bib:slqSwitch:2017} reported that the method proposed by \cite{bib:Xu:2004} took 26 s per upper stage gradient iteration and \cite{bib:slqSwitch:2017} 4.4 s per upper stage gradient iteration on an Intel Core-i7 CPU @2.7 GHz.
As in the previous example, the solution provided by the proposed method is almost identical to the continuous-time counterpart reported in \cite{bib:Xu:2004}; for example, our results show that $t_1 = 0.2335$ [s] and $t_2 = 1.0179$ [s], and the solution of the continuous time was $t_1 = 0.2262$ [s] and $t_2 = 1.0176$ [s].

A future improvement of the proposed method is the Hessian convexification and globalization.
When the initial guess of the switching time is significantly different from the optimal one, the proposed method may diverge or converge to saddle points because the Hessian may contain a negative curvature, that is, Assumption 3.1 may not hold.
For example, the proposed method with only the fraction-to-boundary rule for step-size selection converges to a saddle point when we use the initial guess $t_1 = 1.0$ [s] and $t_2 = 2.0$ [s], for which \cite{bib:Xu:2004} and \cite{bib:slqSwitch:2017} succeed in convergence to the optimal solution with line searches for the upper-stage problem (optimization of the switching instants for a fixed control input).
Furthermore, it is not straightforward to use the existing Hessian convexification method, such as the Gauss-Newton Hessian approximation, because the structure of the Hessian (\ref{eq:QPsub}) is different from that of the standard OCPs.

\section{Conclusion}
We proposed an efficient algorithm for the OCP of nonlinear switched systems that optimizes the control input and switching instants simultaneously for a given switching sequence.
We formulated the OCP based on the direct multiple shooting method with regard to the switching instance as the optimization variable, as well as the state and control input.
We derived a linear equation for Newton's method and a Riccati recursion algorithm to solve the linear equation.
The computational time of the proposed method scales linearly with respect to the length of the horizon as the standard Riccati recursion (\cite{bib:Giaf:2016, bib:Nielsen:2017}).
We conducted numerical experiments and demonstrated that the proposed method converges with a significantly shorter computational time compared with the previous two-stage methods.

As part of our future work, we plan to extend the proposed method to systems with state-dependent switching conditions and state jumps to model mechanical systems that have contact with the environment (\cite{bib:bipedWalking, bib:Katayama:2020, bib:DDP:jumpRobot}).
Subsequently, we need to extend the proposed algorithm to problems with the pure-state constraints that represent the switching conditions, for example, by using the approach proposed by \cite{bib:constrainedRiccati}.
Future work also includes the Hessian convexification and globalization of the proposed method to avoid divergence or convergence to saddle points, even when the initial guess of the solution (particularly the initial guess of the switching instant) is significantly different from the optimal solution.

\bibliography{ifacconf}

\end{document}